\documentclass[reqno]{amsart}
\usepackage[hang]{footmisc}
\usepackage{amssymb,amsthm}
\usepackage{tabularx}
\usepackage{enumerate}
\usepackage{graphicx}
\usepackage{subfigure}
\usepackage{cases}
\usepackage[compress]{cite}
\usepackage{color}
\newtheorem{thm}{Theorem}[section]
\newtheorem{lem}[thm]{Lemma}

\newtheorem{prop}[thm]{Proposition}

\newtheorem{defi}{Definition}[section]
\newtheorem{pppp}{Proof}

\numberwithin{equation}{section}
\newcommand{\bpf}{{\bf Proof:\ \ }}
\newcommand{\epf}{\mbox{}\hfill $\Box$}
\begin{document}
\thispagestyle{empty}

\title{Dynamics of a diffusive competitive model in a periodically evolving domain}

\author{Jiazhen Zhu, Jiazheng Zhou, Zhigui Lin}


\address{Jiazhen Zhu\newline
School of Mathematical Science, Yangzhou University, Yangzhou 225002, China\newline
{\it Email address}: \bf{luckyjiazhenzhu@foxmail.com}}
\address{Jiazheng Zhou\newline
Departamento de Matem\'{a}tica, Universidade de Bras\'{i}lia, BR 70910-900, Bras\'{i}lia-DF, Brazil\newline
{\it Email address}: \bf{zhoumat@hotmail.com}}
\address{Zhigui Lin(corresponding author)\newline
School of Mathematical Science, Yangzhou University, Yangzhou 225002, China\newline
{\it Email address}: \bf{zglin68@hotmail.com; zglin@yzu.edu.cn}}

\maketitle
\medskip
\begin{quote}
\noindent
{Abstract.}{ 
\small In this paper, we are concerned with a two-species competitive model with diffusive terms in a periodically evolving domain and study the impact of the spatial periodic evolution on the dynamics of the model. The Lagrangian transformation approach is adopted to convert the model from a changing domain to a fixed one with the assumption that the evolution of habitat is uniform and isotropic. The ecological reproduction indexes of the linearized model are given as thresholds to reveal the dynamic behaviour of the competitive model. Our theoretical results show that a lager evolving rate benefits the persistence of competitive populations for both sides in the long run. Numerical experiments illustrate that two competitive species, one of which survive and the other vanish in a fixed domain, both survive in a domain with a large evolving rate, and both vanish in a domain with a small evolving rate.}
\end{quote}

\section{Introduction and model formulation}
 A considerable amount of models have been introduced in population ecology. Lotka-Volterra model, a typical population model, was proposed and studied to investigate the behaviour of two species that compete with each other for more survival resources \cite{niwenjie2018}. To understand the possible influence of spatial diffusion which caused by the random movement of individuals within a species, we consider the classic Lotka-Volterra competitive model with diffusive terms $d_{1}\triangle u_1$ and $d_{2}\triangle u_2$ as follows:
\begin{equation}\label{1.1}
\left\{
\begin{array}{lll}
u_{1t}-d_1\triangle u_1=u_1(a_1-c_1u_1-b_1u_2),\ &x\in\Omega,\ t>0,\\
u_{2t}-d_2\triangle u_2=u_2(a_2-b_2u_1-c_2u_2),\ &x\in\Omega,\ t>0,
\end{array} \right.
\end{equation}
where $\Omega\subseteq\mathbb{R}^{n}$ is a non-empty smooth open set, $u_i(x,\ t)(i=1,\ 2)$ represents the density of the $i$-th competitive species depending on location $x$ and time $t$, the positive constant $d_i(i=1,\ 2)$ is the free-diffusion coefficient of $u_i$, and the positive constants $a_i$, $b_i$ and $c_i(i=1,\ 2)$ denote the intrinsic population growth rate, interspecific competition factor and intraspecific competition factor, respectively.

Assume that there is no species across the boundary, the authors in \cite{cantrell2003, pao1992} studied the reaction-diffusive problem:
\begin{equation}
\left\{
\begin{array}{lll}
u_{1t}-d_1\triangle u_1=u_1(a_1-c_1u_1-b_1u_2),\ &x\in\Omega,\ t>0,\\
u_{2t}-d_2\triangle u_2=u_2(a_2-b_2u_1-c_2u_2),\ &x\in\Omega,\ t>0,\\
\frac{\partial u_1(x,\ t)}{\partial \eta}=\frac{\partial u_2(x,\ t)}{\partial \eta}=0,\ &x\in \partial\Omega,\ t>0,\\
u_1(x,\ 0)=u_{1,0}(x),\ u_2(x,\ 0)=u_{2,0}(x),\ &x\in\Omega,
\end{array} \right.
\label{1.2}
\end{equation}
where $\eta$ is the unit outer normal vector of $\partial\Omega$. Clearly, the corresponding steady-state problem of (\ref{1.2}) admits the trivial solution $U_0=(0,\ 0)$ and the semi-trivial solutions $U_1=(\frac{a_1}{c_1},\ 0)$ and $U_2=(0,\ \frac{a_2}{c_2})$. In particular, the steady-state problem admits the unique positive solution $U^*=(\frac{a_1c_2-a_2b_1}{c_1c_2-b_1b_2},\ \frac{a_2c_1-a_1b_2}{c_1c_2-b_1b_2})$ when $\frac{c_1}{b_2}>\frac{a_1}{a_2}>\frac{b_1}{c_2}$ or $\frac{c_1}{b_2}<\frac{a_1}{a_2}<\frac{b_1}{c_2}$. Further theoretical results for stability have been achieved in \cite{pao1992} as follows:

($i$) the trivial solution $U_0=(0,\ 0)$ is always unstable;

$(ii)$ $U^*$ is globally asymptotically stable when $\frac{c_1}{b_2}>\frac{a_1}{a_2}>\frac{b_1}{c_2}$ (weak competition);

$(iii)$ $U_1$ is globally asymptotically stable when $\frac{a_1}{a_2}>\max\{\frac{c_1}{b_2},\ \frac{b_1}{c_2}\}$;

$(iv)$ $U_2$ is globally asymptotically stable when $\frac{a_1}{a_2}<\min\{\frac{c_1}{b_2},\ \frac{b_1}{c_2}\}$;

$(v)$ $U_1$, as well as $U_2$, is locally asymptotically stable and $U^*$ is unstable when $\frac{c_1}{b_2}<\frac{a_1}{a_2}<\frac{b_1}{c_2}$ (strong competition).

Most reaction-diffusion problems describing ecologic models are studied in fixed domains. However, it is common in nature that the habitats in which species live are changeable. Sometimes, boundaries of shifting habitats are unknown owing to the activities of species. For examples, the spreading of invasive species like muskrats in Europe in the early 1900s \cite{Skellam1951}, Asian carps in the Illinois River since the early 1990s \cite{Irons2007}, cane toad (Bufo marinus) in tropical Australia introduced in 1935 \cite{Phillips2007} and the transmission of disease like West Nile virus \cite{Krishnan2016}. Models with such unknown moving boundaries are characterized by free boundary problems and studied as a brunch of model analysis \cite{Friedman}. Mathematically, the free boundary induces more difficulties but it better characterizes the spreading of invasive species\cite{DuLin2010, Lin2014, Lin2015}, and the transformation of disease\cite{DuLin2018, GKLZ, LinZhu2017}. Sometimes, habitat spaces could change following certain known pattern due to objective factors like climate change and seasonal succession. Usually, leaves keep growing before falling and the water storage of lakes annually shifts. For example, the date in \cite{lei2011} give that, in 2009, the wetland vegetation area of Poyang Lake was about 20.8 $km^2$ in February and up to about 1048.9 $km^2$ in May. Fig. 1(a) are the monthly distributions of grassland in Poyang Lake in 2009 from January to December, and Fig. 1(b) is the monthly variation curve of vegetation area \cite{lei2011}. Fig. 1 indicates that the Poyang Lake in China is an evolving domain since the water area of the Lake changes from smaller in winter to larger in summer. Problems with such known boundaries are characterized as growing domain \cite{jorge, madz} or evolving domain \cite{wang2018, xx2016, Sun2020}, and have been studied extensively.

In this paper, we study the Lotka-Volterra competitive model in a periodic evolving domain which refers to a domain evolving with known periodicity.
\begin{figure}
\subfigure[]{
\includegraphics[width=0.4\textwidth]{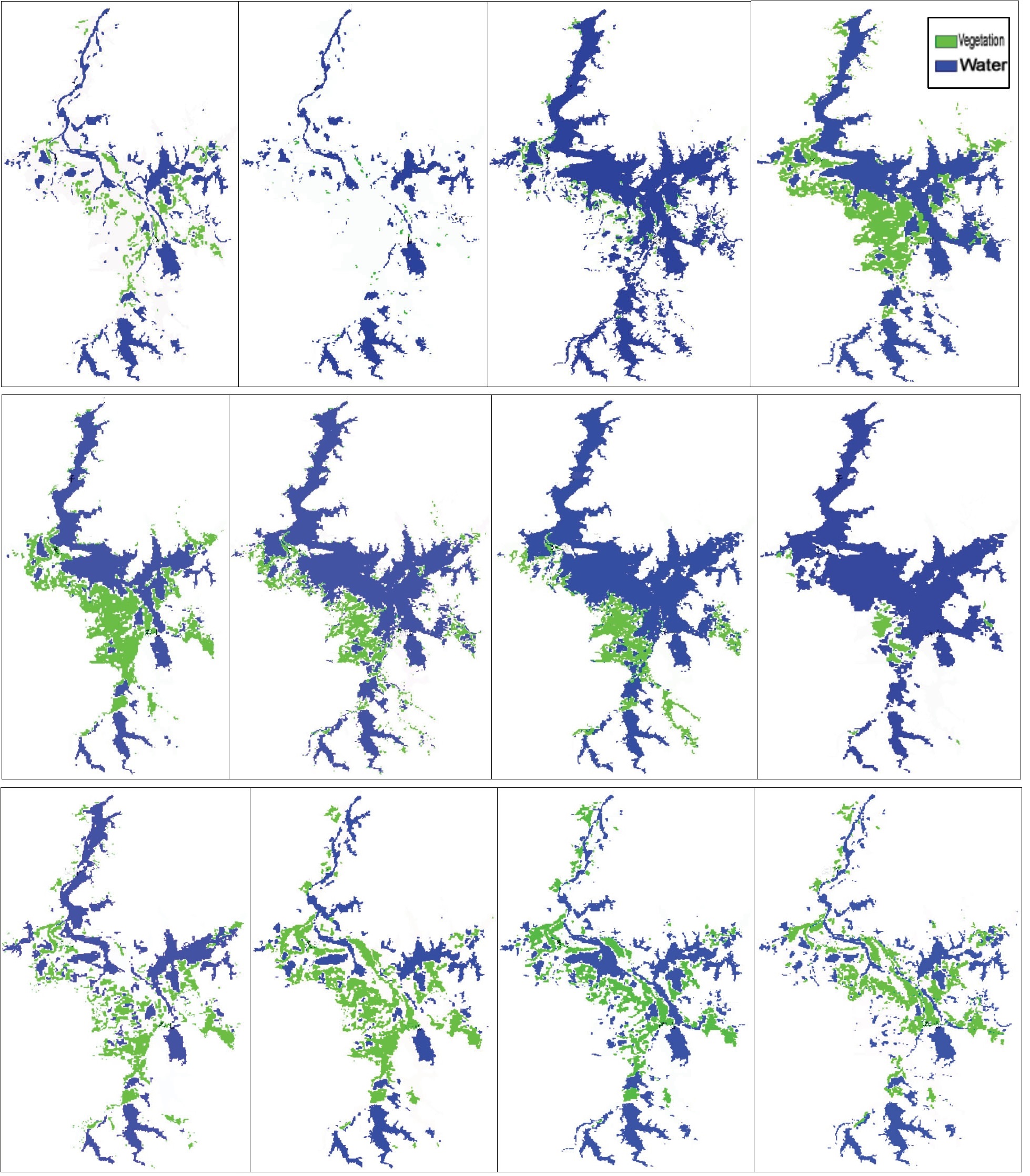}
}
\subfigure[]{
\includegraphics[width=0.56\textwidth]{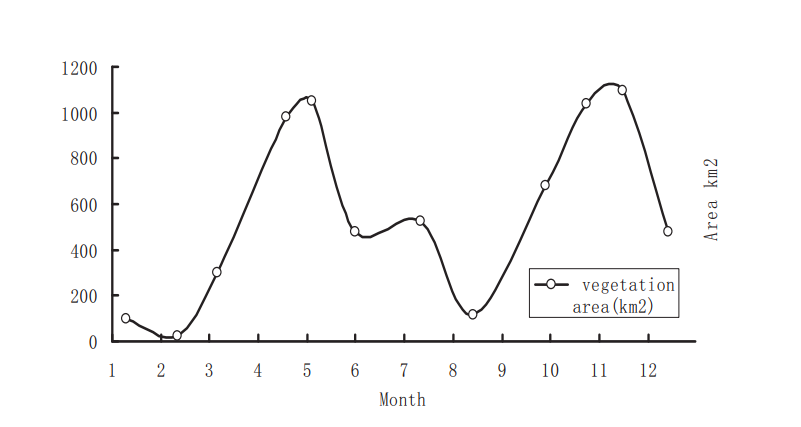}
}
\renewcommand{\figurename}{Fig.}
\caption{\scriptsize (a) are the monthly distribution of grassland and water area in Poyang Lake in 2009 from January to December. (b) is the monthly variation curve of vegetation area which together show the monthly area changes in Poyang Lake\cite{lei2011}.}
\end{figure}
Assume the domain in model (\ref{1.1}) is changing with $t$, that is $\Omega=\Omega(t)\subseteq \mathbb{R}^n$ is time-varying and its boundary $\partial\Omega(t)$ is evolving. According to the principle of mass conservation and Reynolds transport theorem \cite{Acheson1990}, model (\ref{1.1}) can be converted to the following problem in a evolving domain $\Omega(t)$ with Dirichlet boundary condition which implies that there is no species on the boundary:
\begin{equation}
\left\{
\begin{array}{lll}
u_{1t}-d_1\triangle u_1+\textbf{a}\cdot\nabla{u_1}+u_1\nabla\cdot\textbf{a}=u_1(a_1-c_1u_1-b_1u_2),\ &x\in\Omega(t),\ t>0,\\
u_{2t}-d_2\triangle u_2+\textbf{a}\cdot\nabla{u_2}+u_2\nabla\cdot\textbf{a}=u_2(a_2-b_2u_1-c_2u_2),\ &x\in\Omega(t),\ t>0,\\
u_1(x(t),\ t)=u_2(x(t),\ t)=0,\ &x\in\partial\Omega(t),\ t>0,\\
u_1(x(0),\ 0)=u_{1,0}(x(0)),\ u_2(x(0),\ 0)=u_{2,0}(x(0)),\ &x(0)\in\overline\Omega(0),
\end{array} \right.
\label{1.3}
\end{equation}
where \textbf{a} denotes the spacial flow velocity caused by the change of domain, $u_1\cdot\nabla\textbf{a}$ and $u_2\cdot\nabla\textbf{a}$ are called dilution terms, $\textbf{a}\cdot\nabla{u_1}$ and $\textbf{a}\cdot\nabla{u_2}$ are called advection terms. $x=x(t)$ within $\overline{\Omega}(t)$ is the function of $t$, $a_i=a_i(t),\ b_i=b_i(t)$ and $c_i=c_i(t)(i=1,\ 2)$ are all positive and $T$-periodic.

Assume the evolution of $\Omega(t)$ is uniform and isotropic, that is,
\begin{equation}\label{l.4}
   x(t)=\rho(t)y,\ y\in\Omega(0),
\end{equation}
where $\rho(t)=\rho(t+T)$ is a $T$-periodic function with $\rho(0)=1$.
Thus, $u_1$ and $u_2$ can be mapped as a new function with the definition:
\begin{equation}\label{1.5}
   u_1(x,\ t)=v_1(y,\ t),\ u_2(x,\ t)=v_2(y,\ t)
\end{equation}
followed with
\begin{equation*}
\left.
\begin{array}{lll}
   v_{1t}=\frac{\partial u_1}{\partial t}+\textbf{a}\cdot\nabla{u_1},\ v_{2t}=\frac{\partial u_2}{\partial    t}+\textbf{a}\cdot\nabla{u_2},\\
   \nabla\textbf{a}=\frac{n\dot{\rho}(t)}{\rho(t)},\\
   \Delta u_1=\frac{1}{\rho^2(t)}\Delta v_1,\ \Delta u_2=\frac{1}{\rho^2(t)}\Delta v_2,
\end{array}
\right.
\end{equation*}
where $n$ is the dimension of the space $\Omega$.

Therefore, (\ref{1.3}) is converted to the problem in a fixed domain
\begin{equation}
\left\{
\begin{array}{lll}
v_{1t}-\frac{d_1}{\rho^2(t)}\triangle v_1=-\frac{n\dot{\rho}(t)}{\rho(t)}v_1+v_1(a_1-c_1v_1-b_1v_2),\ &y\in\Omega(0),\ t>0,\\
v_{2t}-\frac{d_2}{\rho^2(t)}\triangle v_2=-\frac{n\dot{\rho}(t)}{\rho(t)}v_2+v_2(a_2-b_2v_1-c_2v_2),\ &y\in\Omega(0),\ t>0,\\
v_1(y,\ t)=v_2(y,\ t)=0,\ &y\in\partial\Omega(0),\ t>0,\\
v_1(y,\ 0)=v_{1,0}(y),\ v_2(y,\ 0)=v_{2,0}(y),\ &y\in\overline\Omega(0),
\end{array} \right.
\label{1.6}
\end{equation}
the dynamics of which is related to its corresponding periodic problem
\begin{equation}
\left\{
\begin{array}{lll}
V_{1t}-\frac{d_1}{\rho^2(t)}\triangle V_1=-\frac{n\dot{\rho}(t)}{\rho(t)}V_1+V_1(a_1-c_1V_1-b_1V_2),\ &y\in\Omega(0),\ t>0,\\
V_{2t}-\frac{d_2}{\rho^2(t)}\triangle V_2=-\frac{n\dot{\rho}(t)}{\rho(t)}V_2+V_2(a_2-b_2V_1-c_2V_2),\ &y\in\Omega(0),\ t>0,\\
V_1(y,\ t)=V_2(y,\ t)=0,\ &y\in\partial\Omega(0),\ t>0,\\
V_1(y,\ 0)=V_1(y,\ T),\ V_2(y,\ 0)=V_2(y,\ T),\ &y\in\overline\Omega(0).
\end{array} \right.
\label{1.7}
\end{equation}

In the rest of this paper, we are devoted to investigating the asymptotic behaviour of the initial and boundary value problem (\ref{1.6}) in related to the $T$-periodic solution of problem (\ref{1.7}). In Section 2, we first present the ecological reproduction indexes of problem (\ref{1.7}) as thresholds based on the principal eigenvalues of its linearized problem, and then deliver the existence of periodic solution. In Section 3, we analyze the stability of the solution to the initial and boundary value problem. In Section 4, we discuss the impact of the evolving domain on the persistence of two competitive species. In Section 5, we give some numerical simulations and ecological explanations in support of the theoretical results achieved in Section 4.

\section{Ecological reproduction index}
In this section, we are going to determine the existence of the solution to problem (\ref{1.7}). After linearizing problem (\ref{1.7}) around $(0,\ 0)$, we have its eigenvalue problem as follows:
\begin{equation}
\left\{
\begin{array}{lll}
\phi_{1t}-\frac{d_1}{\rho^2(t)}\triangle \phi_1=(a_1-\frac{n\dot{\rho}(t)}{\rho(t)})\phi_1+\lambda_1\phi_1,
\ &y\in\Omega(0),\ t>0,\\
\phi_{2t}-\frac{d_2}{\rho^2(t)}\triangle \phi_2=(a_2-\frac{n\dot{\rho}(t)}{\rho(t)})\phi_2+\lambda_2\phi_2,
\ &y\in\Omega(0),\ t>0,\\
\phi_1(y,\ t)=\phi_2(y,\ t)=0,\ &y\in\partial\Omega(0),\ t>0,\\
\phi_1(y,\ 0)=\phi_1(y,\ T),\ \phi_2(y,\ 0)=\phi_2(y,\ T),\ &y\in\overline\Omega(0),
\end{array} \right.
\label{2.1}
\end{equation}
and denote $\lambda^\vartriangle_i(i=1,\ 2)$ the principle eigenvalue of (\ref{2.1}), and $\phi^\vartriangle_i$ the corresponding eigenfunctions with $0\leq \phi^\vartriangle_i\leq 1$. Furthermore, by variation method, we can give the explicit expression of principal eigenvalues as:
\begin{center}
 ${\lambda}^\vartriangle_i =\frac{1}{T}\int_{0}^{T}{\frac{d_{i}\lambda_{0}}{\rho^{2}(t)}}dt-\frac{1}{T}\int_{0}^{T}{a_{i}(t)}dt\ (i=1,\ 2)$,
\end{center}
where $\lambda_0$ is the principal eigenvalue to
\begin{equation}
\left\{
\begin{array}{lll}
\ -\triangle \phi=\lambda \phi,
\ &y\in\Omega(0),\\
\phi(y)=0,\ &y\in\partial\Omega(0).\\
\end{array} \right.
\label{2.2}
\end{equation}
Using the next generation operator as in \cite{LiangX, Zhaobook}, we can define the ecological reproduction index $R_i(i=1,\ 2)$. Moreover,
it follows from Lemma 13.1.1 in  \cite{Zhaobook} that $R_1$ and $R_2$ are the principal eigenvalues of the following problems:
\begin{equation}
\left\{
\begin{array}{lll}
\varphi_{1t}-\frac{d_1}{\rho^2(t)}\triangle \varphi_1=(\frac{a_1}{R_1}-\frac{n\dot{\rho}(t)}{\rho(t)})\varphi_1,
\ &y\in\Omega(0),\ t>0,\\
\varphi_{2t}-\frac{d_2}{\rho^2(t)}\triangle \varphi_2=(\frac{a_2}{R_2}-\frac{n\dot{\rho}(t)}{\rho(t)})\varphi_2,
\ &y\in\Omega(0),\ t>0,\\
\varphi_1(y,\ t)=\varphi_2(y,\ t)=0,\ &y\in\partial\Omega(0),\ t>0,\\
\varphi_1(y,\ 0)=\phi_1(y,\ T),\ \varphi_2(y,\ 0)=\phi_2(y,\ T),\ &y\in\overline\Omega(0).
\end{array} \right.
\label{2.3}
\end{equation}
In the study of epidemic model, $R_i$ is called basic reproduction number \cite{Dietz1993, LiangX} and usually given as threshold. Similarly, the variation method gives
\begin{equation}\label{2.4}
  R_i=\frac{\int_{0}^{T}a_i(t)dt}{d_i\lambda_0\int_{0}^{T}\frac{1}{\rho^2(t)}dt}\, (i=1,\ 2).
\end{equation}
It can be verified that
\begin{equation}
   {\rm sgn}(1-R_i)={\rm sgn}(\lambda_i)\, (i=1,\ 2).
\label{2.5}
\end{equation}
Similar results for general systems hold as well. For more details, see \cite{LiangX} and the references therein.
To derive the existence of the solution to (\ref{1.7}), we give the definition of upper and lower solutions.
\begin{defi}\label{def2.1}
$(\tilde{V}_1,\ \tilde{V}_2)$ and $(\hat{V}_1,\ \hat{V}_2)$ is a pair of coupled upper and lower solutions of the problem $(\ref{1.7})$, if
\begin{equation}
\left\{
\begin{array}{lll}
\hat{V}_{1t}-\frac{d_1}{\rho^2(t)}\triangle \hat{V}_1\leq-\frac{n\dot{\rho}(t)}{\rho(t)}\hat{V}_1+\hat{V}_1(a_1-c_1\hat{V}_1-b_1\tilde{V}_2),\ &y\in\Omega(0),\ t>0,\\
\tilde{V}_{1t}-\frac{d_1}{\rho^2(t)}\triangle \tilde{V}_1\geq-\frac{n\dot{\rho}(t)}{\rho(t)}\tilde{V}_1+\tilde{V}_1(a_1-c_1\tilde{V}_1-b_1\hat{V}_2),\ &y\in\Omega(0),\ t>0,\\
\hat{V}_{2t}-\frac{d_2}{\rho^2(t)}\triangle \hat{V}_2\leq-\frac{n\dot{\rho}(t)}{\rho(t)}\hat{V}_2+\hat{V}_2(a_2-b_2\tilde{V}_1-c_2\hat{V}_2),\ &y\in\Omega(0),\ t>0,\\
\tilde{V}_{2t}-\frac{d_2}{\rho^2(t)}\triangle \tilde{V}_2\geq-\frac{n\dot{\rho}(t)}{\rho(t)}\tilde{V}_2+\tilde{V}_2(a_2-b_2\hat{V}_1-c_2\tilde{V}_2),\ &y\in\Omega(0),\ t>0,\\
\tilde{V}_1(y,\ t)\geq \hat{V}_1(y,\ t)=0,\ \tilde{V}_2(y,\ t)\geq \hat{V}_2(y,\ t)=0,\ &y\in\partial\Omega(0),\ t\geq0,\\
\hat{V}_1(y,\ 0)\leq \hat{V}_1(y,\ T),\ \hat{V}_2(y,\ 0)\leq \hat{V}_2(y,\ T),\ &y\in\overline\Omega(0),\\
\tilde{V}_1(y,\ 0)\geq \tilde{V}_1(y,\ T),\ \tilde{V}_2(y,\ 0)\geq \tilde{V}_2(y,\ T),\ &y\in\overline\Omega(0).
\end{array} \right.
\end{equation}
\end{defi}
Let $S_0:=\{\,(V_1,\ V_2):\ (\hat{V}_1,\ \hat{V}_2)\leq(V_1,\ V_2)\leq(\tilde{V}_1,\ \tilde{V}_2),\ (y,\, t)\in\overline{\Omega}(0)\times[0,\ T]\,\}$ and denote
\begin{center}
  $f_1(V_1,\ V_2)=V_1(a_1-c_1V_1-b_1V_2)-\frac{n\dot{\rho(t)}}{\rho(t)}V_1,$
\end{center}
and
\begin{center}
  $f_2(V_1,\ V_2)=V_2(a_2-b_2V_1-c_2V_2)-\frac{n\dot{\rho(t)}}{\rho(t)}V_2.$
\end{center}
Then, for any $(V_1,\ V_2),\ (Z_1,\ Z_2)\in S_0$,
{\small
\begin{equation*}
\left|f_{1}(V_{1},\ V_{2})-f_{1}(Z_1,\ Z_2)\right|
\leq[a_{1}^{M}+(b_{1}^{M}+2c_{1}^{M})\frac{a_{1}^{M}}{c_{1}^{m}}+b_{1}^{M}\frac{a_{2}^{M}}{c_{2}^{m}}+\frac{n\dot{\rho}^{M}}{\rho^{m}}]
(\left|V_{1}-Z_{1}\right|+\left|V_{2}-Z_{2}\right|),
\end{equation*}
\begin{equation*}
\left|f_{2}(V_{1},\ V_{2})-f_{2}(Z_{1},\ Z_{2})\right|
\leq[a_{2}^{M}+(b_{2}^{M}+2c_{2}^{M})\frac{a_{2}^{M}}{c_{2}^{m}}+b_{2}^{M}\frac{a_{1}^{M}}{c_{1}^{m}}+\frac{n\dot{\rho}^{M}}{\rho^{m}}]
(\left|V_{1}-Z_{1}\right|+\left|V_{2}-Z_{2}\right|),
\end{equation*}}
where $f^M=\max\limits_{[0,\ T]}f(t)$ and $f^m=\min\limits_{[0,\ T]}f(t)$. We find that $f_1$ and $f_2$ satisfy the Lipschitz condition with Lipschitz coefficients
\begin{equation}
  k_{1}=a_{1}^{M}+(b_{1}^{M}+2c_{1}^{M})\frac{a_{1}^{M}}{c_{1}^{m}}+b_{1}^{M}\frac{a_{2}^{M}}{c_{2}^{m}}+\frac{n\dot{\rho}^{M}}{\rho^{m}},
\end{equation}
and
\begin{equation}
  k_{2}=a_{2}^{M}+(b_{2}^{M}+2c_{2}^{M})\frac{a_{2}^{M}}{c_{2}^{m}}+b_{2}^{M}\frac{a_{1}^{M}}{c_{1}^{m}}+\frac{n\dot{\rho}^{M}}{\rho^{m}}.
\end{equation}
Based on the upper and lower solutions technique developed by Pao \cite{Pao2005}, we have the following result about the existence of the solution.
\begin{lem}\label{lem2.1}
If $(\tilde{V}_1,\ \tilde{V}_2),\ (\hat{V}_1,\ \hat{V}_2)$ is a pair of coupled upper and lower solutions of $(\ref{1.7})$, then $(\ref{1.7})$ admits at least one periodic solution $(V_1,\ V_2)\in S_0$.
\end{lem}

Now we present the existence of the periodic solution to (\ref{1.7}).
\begin{thm}\label{thm2.1}
Denote $M_1=(\frac{1}{c_1}(a_1-\frac{n\dot{\rho}}{\rho}))^M$ and $M_2=(\frac{1}{c_2}(a_2-\frac{n\dot{\rho}}{\rho}))^M$. Then we have the following assertions:

$(i)$ if $R_1\leq1$ and $R_2\leq1$, $(1.7)$ admits only trivial solution $(0,\ 0)$;

$(ii)$ if $R_1>1$ and $R_2\leq1$, $(1.7)$ admits a semi-trivial periodic solution $(V_1^{\vartriangle},\ 0)$;

$(iii)$ if $R_1\leq1$ and $R_2>1$, $(1.7)$ admits a semi-trivial periodic solution $(0,\ V_2^{\vartriangle})$;

$(iv)$ if $R_1>1$ and $R_2>1$, together with $(\frac{a_1}{b_1})^m(1-\frac{1}{R_1})>M_2$ and $(\frac{a_2}{b_2})^m(1-\frac{1}{R_2})>M_1$,  $(1.7)$ admits a positive periodic solution $(V_1^{*},\ V_2^*)$.
\end{thm}
\bpf
$(i)$ Let $(V_1,\ V_2)$ be the nonnegative solution of (\ref{1.7}), we claim that $V_1\equiv0$ and $V_2\equiv0$ in $\overline{\Omega}(0)$. In fact, assume that $V_1$ satisfies
\begin{center}
  $V_{1t}-\frac{d_1}{\rho^2(t)}\triangle V_1+\frac{n\dot{\rho}(t)}{\rho(t)}V_1-a_1V_1=-(b_1V_2+c_1V_1)V_1,\ y\in\Omega(0),\ t>0,$
\end{center}
and $V_1\geq0(\not\equiv0)$ by contradiction.
Recalling that
\begin{center}
  $\phi_{1t}-\frac{d_1}{\rho^2(t)}\triangle \phi_1+\frac{n\dot{\rho}(t)}{\rho(t)}\phi_1-a_1\phi_1=\lambda_1^{\vartriangle}\phi_1,\ y\in\Omega(0),\ t>0,$
\end{center}
we have $\lambda_1^{\vartriangle}<0$ according to the monotonicity of eigenvalues revealed in \cite{pengrui}(Proposition 5.2). It follows from (\ref{2.3}) that $\lambda_1^{\vartriangle}<0$ implies $R_1>1$, which leads a contradiction to the condition. Therefore, $V_1\equiv0$ in $\overline{\Omega}(0)$. Similarly, $V_2\equiv0$. Thus, $(0,\ 0)$ is the only nonnegative solution to (1.7).

$(ii)$ If $R_1>1$ and $R_2\leq1$, consider semi-trivial solution $(V_1^{\vartriangle},\ 0)$  and $V_1^{\vartriangle}$ satisfies
\begin{equation}
\left\{
\begin{array}{lll}
  V_{1,\ t}^{\vartriangle}-\frac{d_1}{\rho^2(t)}\Delta V_1^{\vartriangle}=-\frac{n\dot{\rho}(t)}{\rho(t)}V_1^{\vartriangle}+V_1^{\vartriangle}(a_1-c_1V_1^{\vartriangle})\ &y\in\Omega(0),\ t>0,\\
  V_1^{\vartriangle}=0\ &y\in\partial\Omega(0),\ t\geq0,\\
  V_1^{\vartriangle}(y,\ 0)=V_1^{\vartriangle}(y,\ T)\ &y\in\overline\Omega(0).
\end{array} \right.
\label{2.6}
\end{equation}
It can be verified that $M_1$ and $\delta\varphi_1$ is a pair of ordered upper and lower solutions of problem (\ref{2.6}) for any positive constant $\delta< -\lambda_1^{\vartriangle}$. Furthermore, according to Theorem 27.1 in \cite{hess} for the uniqueness of the solution to a problem with concave nonlinearities, the positive solution $V_1^{\vartriangle}$ is unique as $a_1-c_1V_1^{\vartriangle}$ is monotone decreasing in terms of $V_1^{\vartriangle}$. Thus, $(V_1^{\vartriangle},\ 0)$ is the unique periodic solution of (\ref{1.7}).

$(iii)$ The proof of $(iii)$ is similar to that of $(ii)$.

$(iv)$ According to Lemma \ref{lem2.1}, (\ref{1.7}) admits at least one periodic solution $(V_1,\ V_2)$ if we can verify that $(M_1,\ M_2)$ and $(\varepsilon\varphi_1,\ \varepsilon\varphi_2)$
is a pair of coupled upper and lower solutions of (\ref{1.7}) with positive constant $\varepsilon$ to be determined.
In fact, the choose of $M_1$ and $M_2$ implies that $(M_1,\ M_2)$ is an upper solution of (\ref{1.7}) as long as $(\varepsilon\varphi_1,\ \varepsilon\varphi_2)$ is nonnegative. Clearly, the condition $(\frac{a_1}{b_1})^m(1-\frac{1}{R_1})>M_2$ and $(\frac{a_2}{b_2})^m(1-\frac{1}{R_2})>M_1$ implies that there exists a constant
\begin{center}
 $\varepsilon_0={\rm min}\{\frac{1}{c_1^M}(a_1(1-\frac{1}{R_1})-b_1M_2),\ \frac{1}{c_2^M}(a_2(1-\frac{1}{R_2})-b_2M_1)\}>0$,
\end{center}
then for any $0<\varepsilon<\varepsilon_0$, $(\varepsilon\varphi_1,\ \varepsilon\varphi_2)$ is the lower solution of (\ref{1.7}) with $(M_1,\ M_2)$ the upper solution. Thus, $(M_1,\ M_2)$ and $(\varepsilon\varphi_1,\ \varepsilon\varphi_2)$ is a pair of coupled upper and lower solutions of (\ref{1.7}) and the proof is completed.
\epf

\section{Dynamics of periodic solutions}
In this section, we are going to discuss the stability of the solution to problem (\ref{1.6}) which is related to the solution of the periodic problem (\ref{1.7}). Firstly, we convert the reaction functions in problem (\ref{1.6}) to be quasimonotone nondecreasing.

Let $M=\max\{M_2,\ \sup\limits_{y\in\Omega(0)}v_{2,0}(y)\}$, $v_3=M-v_2$ and then (\ref{1.6}) becomes
\begin{equation}
\left\{
\begin{array}{lll}
v_{1t}-\frac{d_1}{\rho^2(t)}\Delta v_1=f_1(v_1,\ M-v_3),\ &y\in\Omega(0),\ t>0,\\
v_{3t}-\frac{d_2}{\rho^2(t)}\Delta v_3=-f_2(v_1,\  M-v_3),\ &y\in\Omega(0),\ t>0,\\
v_1(y,\ t)=M-v_3(y,\ t)=0,\ &y\in\partial\Omega(0),\ t>0,\\
v_1(y,\ 0)=v_{1,0}(y),\ &y\in\overline\Omega(0),\\
v_3(y,\ 0)=v_{3,0}(y):=M-v_{2,0}(y),\ &y\in\overline\Omega(0).
\end{array} \right.
\label{3.1}
\end{equation}
The corresponding periodic problem of (\ref{3.1}) becomes
\begin{equation}
\left\{
\begin{array}{lll}
V_{1t}-\frac{d_1}{\rho^2(t)}\triangle V_1=f_1(V_1,\ M-V_3),\ &y\in\Omega(0),\ t>0,\\
V_{3t}-\frac{d_2}{\rho^2(t)}\triangle V_3=-f_2(V_1,\ M-V_3),\ &y\in\Omega(0),\ t>0,\\
V_1(y,\ t)=M-V_3(y,\ t)=0,\ &y\in\partial\Omega(0),\ t>0,\\
V_1(y,\ 0)=V_1(y,\ T),\ V_3(y,\ 0)=V_3(y,\ T),\ &y\in\overline\Omega(0),
\end{array} \right.
\label{3.2}
\end{equation}
where
\begin{center}
  $f_1(V_1,\ M-V_3)=-\frac{n\dot{\rho}(t)}{\rho(t)}V_1+V_1(a_1-b_1(M-V_3)-c_1V_1)$
\end{center}
and
\begin{center}
  $-f_2(V_1,\ M-V_3)=\frac{n\dot{\rho}(t)}{\rho(t)}(M-V_3)-(M-V_3)(a_2-b_2V_1-c_2(M-V_3)))$
\end{center}
are quasimonotone nondecreasing reaction functions for $(V_1,\ V_3)\in S_1$, where
\begin{center}
  $S_1:=\{(V,\ Z):\ (\hat{V}_1,\ M-\tilde{V}_2)\leq(V,\ Z)\leq(\tilde{V}_1,\ M-\hat{V}_2),\ (y, t)\in\overline{\Omega}(0)\times[0,\ T]\}$.
\end{center}
We claim that $(\tilde{V_1},\ M-\hat{V_2})$ and $(\hat{V_1},\ M-\tilde{V_2})$ is a pair of ordered upper and lower solutions of (\ref{3.2}) if $(\tilde{V_1},\ \tilde{V_2})$ and $(\hat{V_1},\ \hat{V_2})$ is a pair of coupled nonnegative upper and lower solutions of (\ref{1.7}). And sequences $\{(\overline{V}_1^{(m)},\ \overline{V}_3^{(m)})\}$ and $\{(\underline{V}_1^{(m)},\ \underline{V}_3^{(m)})\}$ can be obtained by taking $\overline{V}_1^{(0)}=\tilde{V_1},\ \overline{V}_3^{(0)}=M-\hat{V_2},\ \underline{V}_1^{(0)}=\hat{V_1}$ and $\underline{V}_3^{(0)}=M-\tilde{V_2}$ as initial iterations and solving the linear periodic problem
\begin{equation}
\left\{
\begin{array}{lll}
\overline{V}_{1t}^{(n)}-\frac{d_1}{\rho^2(t)}\Delta \overline{V}_1^{(n)}+k_1\overline{V}_1^{(n)}=F_1(t,\ \overline{V}_1^{(n-1)},\ \overline{V}_3^{(n-1)}),\ &y\in\Omega(0),\ t>0,\\
\overline{V}_{3t}^{(n)}-\frac{d_2}{\rho^2(t)}\Delta \overline{V}_3^{(n)}+k_2\overline{V}_3^{(n)}=F_2(t,\ \overline{V}_1^{(n-1)},\ \overline{V}_3^{(n-1)}),\ &y\in\Omega(0),\ t>0,\\
\underline{V}_{1t}^{(n)}-\frac{d_1}{\rho^2(t)}\Delta \underline{V}_1^{(n)}+k_1\underline{V}_1^{(n)}=F_1(t,\ \underline{V}_1^{(n-1)},\ \underline{V}_3^{(n-1)}),\ &y\in\Omega(0),\ t>0,\\
\underline{V}_{3t}^{(n)}-\frac{d_2}{\rho^2(t)}\Delta \underline{V}_3^{(n)}+k_2\underline{V}_3^{(n)}=F_2(t,\ \underline{V}_1^{(n-1)},\ \underline{V}_3^{(n-1)}),\ &y\in\Omega(0),\ t>0,\\
\overline{V}_1^{(n)}=\underline{V}_1^{(n)}=0,\ \overline{V}_3^{(n)}=\underline{V}_3^{(n)}=M,\ &y\in\partial\Omega(0),\ t>0,\\
\overline{V}_1^{(n)}(y,\ 0)=\overline{V}_1^{(n-1)}(y,\ T),\ \overline{V}_3^{(n)}(y,\ 0)=\overline{V}_3^{(n-1)}(y,\ T),\ &y\in\overline\Omega(0),\\
\underline{V}_1^{(n)}(y,\ 0)=\underline{V}_1^{(n-1)}(y,\ T),\ \underline{V}_3^{(n)}(y,\ 0)=\underline{V}_3^{(n-1)}(y,\ T),\ &y\in\overline\Omega(0),
\end{array}
\right.
\label{3.3}
\end{equation}
where $k_1$ and $k_2$ are Lipschitz coefficients given in (2.7) and (2.8),
\begin{center}
  $F_1(t,\ V_1,\ V_3)=k_1V_1+f_1(t,\ V_1,\ M-V_3)$
\end{center}
and
\begin{center}
  $F_2(t,\ V_1,\ V_3)=k_2V_3-f_2(t,\ V_1,\ M-V_3)$.
\end{center}

Similarly, the sequences $\{(\overline{v}_1^{(m)},\ \overline{v}_3^{(m)})\}$ and $\{(\underline{v}_1^{(m)},\ \underline{v}_3^{(m)})\}$ can be obtained by taking $\overline{v}_1^{(0)}=\tilde{v}_1,\ \overline{v}_3^{(0)}=M-\hat{v}_2,\ \underline{v}_1^{(0)}=\hat{v}_1$ and $\underline{v}_3^{(0)}=M-\tilde{v}_2$ as initial iterations and solving the linear initial and boundary value problem
\begin{equation}
\left\{
\begin{array}{lll}
\overline{v}_{1t}^{(n)}-\frac{d_1}{\rho^2(t)}\Delta \overline{v}_1^{(n)}+k_1\overline{v}_1^{(n)}=F_1(t,\ \overline{v}_1^{(n-1)},\ \overline{v}_3^{(n-1)}),\ &y\in\Omega(0),\ t>0,\\
\overline{v}_{3t}^{(n)}-\frac{d_2}{\rho^2(t)}\Delta \overline{v}_3^{(n)}+k_2\overline{v}_3^{(n)}=F_2(t,\ \overline{v}_1^{(n-1)},\ \overline{v}_3^{(n-1)}),\ &y\in\Omega(0),\ t>0,\\
\underline{v}_{1t}^{(n)}-\frac{d_1}{\rho^2(t)}\Delta \underline{v}_1^{(n)}+k_1\underline{v}_1^{(n)}=F_1(t,\ \underline{v}_1^{(n-1)},\ \underline{v}_3^{(n-1)}),\ &y\in\Omega(0),\ t>0,\\
\underline{v}_{3t}^{(n)}-\frac{d_2}{\rho^2(t)}\Delta \underline{v}_3^{(n)}+k_2\underline{v}_3^{(n)}=F_2(t,\ \underline{v}_1^{(n-1)},\ \underline{v}_3^{(n-1)}),\ &y\in\Omega(0),\ t>0,\\
\overline{v}_1^{(n)}=\underline{v}_1^{(n)}=0,\ \overline{v}_3^{(n)}=\underline{v}_3^{(n)}=M,\ &y\in\partial\Omega(0),\ t>0,\\
\overline{v}_1^{(n)}(y,\ 0)=\underline{v}_1^{(n)}(y,\ 0)=v_{1,0}(y),\ &y\in\overline\Omega(0),\\
\overline{v}_3^{(n)}(y,\ 0)=\underline{v}_3^{(n)}(y,\ 0)=v_{3,0}(y),\ &y\in\overline\Omega(0),
\end{array}
\right.
\label{3.4}
\end{equation}
where $(v_{1,0}(y),\ v_{3,0}(y))\in S_1$.

Next, we present two propositions about the sequences
\begin{center}
$\{(\overline{V}_1^{(m)},\ \overline{V}_3^{(m)})\}$, $\{(\underline{V}_1^{(m)},\ \underline{V}_3^{(m)})\}$, $\{(\overline{v}_1^{(m)},\ \overline{v}_3^{(m)})\}$ and $\{(\underline{v}_1^{(m)},\ \underline{v}_3^{(m)})\}$
\end{center}
according to Pao's work in \cite{Pao2005}.
\begin{prop}\label{prop3.1}
$(i)$ The sequence $\{(\overline{V}_1^{(m)},\ \overline{V}_3^{(m)})\}$ decreases and converges monotonically to $(\overline{V}_1,\ \overline{V}_3)$ which is a maximal $T$-periodic solution of $(\ref{3.2})$, and the sequence $\{(\underline{V}_1^{(m)},\ \underline{V}_3^{(m)})\}$ increases and converges monotonically to $(\underline{V}_1,\ \underline{V}_3)$ which is a minimal $T$-periodic solution of $(\ref{3.2})$, that is
\begin{center}
  $(\hat{V}_1,\ M-\tilde{V}_2)\leq(\underline{V}_1^{(m)},\ \underline{V}_3^{(m)})\leq(\underline{V}_1^{(m+1)},\ \underline{V}_3^{(m+1)})\leq(\underline{V}_1,\ \underline{V}_3)$\\
  $\leq(\overline{V}_1,\ \overline{V}_3)\leq(\overline{V}_1^{(m+1)}, \overline{V}_3^{(m+1)})\leq(\overline{V}_1^{(m)},\ \overline{V}_3^{(m)})\leq(\tilde{V}_1,\ M-\hat{V}_2)$.
\end{center}

$(ii)$ $(\underline{V}_1,\ \underline{V}_3)=(\overline{V}_1,\ \overline{V}_3)$ when $\overline{V}_1(y,\ 0)=\underline{V}_1(y,\ 0)$ and $\overline{V}_3(y,\ 0)=\underline{V}_3(y,\ 0)$ which implies that $(\ref{3.2})$ admits a unique periodic solution
\begin{center}
  $({V}_1,\ {V}_3)(=(\overline{V}_1,\ \overline{V}_3)=(\underline{V}_1,\ \underline{V}_3))$.
\end{center}
\end{prop}
\begin{prop}\label{prop3.2}
Both $\{(\overline{v}_1^{(m)},\ \overline{v}_3^{(m)})\}$ and $\{(\underline{v}_1^{(m)},\ \underline{v}_3^{(m)})\}$ converge to $(v_1,\ v_3)$, the unique solution of $(\ref{3.1})$ satisfying
\begin{center}
   $(\hat{V}_1,\ M-\tilde{V}_2)\leq(\underline{v}_1^{(m)},\ \underline{v}_3^{(m)})\leq(\underline{v}_1^{(m+1)},\ \underline{v}_3^{(m+1)})\leq(v_1,\ v_3)$\\
   $\leq(\overline{v}_1^{(m+1)},\ \overline{v}_3^{(m+1)})\leq(\overline{v}_1^{(m)},\ \overline{v}_3^{(m)})\leq(\tilde{V}_1,\ M-\hat{V}_2)$.
\end{center}
\end{prop}

Based on Propositions \ref{prop3.1} and \ref{prop3.2}, we have the following lemma and detailed proof for more general parabolic systems can be found in \cite{Pao2005}.
\begin{lem}\label{prop3.3}
Let $\boldsymbol{\eta}=(v_{1,0}(y),\ v_{3,0}(y))$ and for any $m$ and $m^{'}$, if
\begin{center}
$(\underline{V}_1^{(m^{'})},\ \underline{V}_3^{(m^{'})})(y,0)\leq\boldsymbol{\eta}(y)\leq(\overline{V}_1^{(m)},\ \overline{V}_3^{(m)})(y,0)$,
\end{center}
then we have that

$(i)$ $(\overline{V}_1^{(m)},\ \overline{V}_3^{(m)})$ and $(\underline{V}_1^{(m^{'})},\ \underline{V}_3^{(m^{'})})$ is a pair of ordered upper and lower solutions of problem $(\ref{3.1})$;

$(ii)$ the solution of $(\ref{3.1})$ denoted by $\mathbf{v}(y,\ t;\ \boldsymbol{\eta})$ satisfies
\begin{center}
  $(\underline{V}_1^{(m)},\ \underline{V}_3^{(m)})(y,\ t)\leq\mathbf{v}(y,\ t+mT;\ \boldsymbol{\eta})\leq(\overline{V}_1^{(m)},\ \overline{V}_3^{(m)})(y,\ t)$
\end{center}
with
\begin{equation}\label{3.5}
\begin{array}{lll}
(\underline{V}_1,\ \underline{V}_3)(y,\ t)&\leq \liminf\limits_{m\rightarrow+\infty} \mathbf{v}(y,\ t+mT;\ \boldsymbol{\eta})\\
&\leq \limsup\limits_{m\rightarrow+\infty}\mathbf{v}(y,\ t+mT;\ \boldsymbol{\eta})\leq (\overline{V}_1,\ \overline{V}_3)(y,\ t).
\end{array}
\end{equation}
\end{lem}

\begin{thm}\label{thm3.4}
Denote $V_3^{\vartriangle}=M-V_2^{\vartriangle}$.
For problem \eqref{3.1} with any nonnegative nontrivial initial value $\boldsymbol{\eta}$, we have the following stability results:

$(i)$ If $R_1\leq1$ and $R_2\leq1$, then $\lim\limits_{m\rightarrow\infty}\mathbf{v}(y,\ t+mT;\ \boldsymbol{\eta})=(0,\ M)$;

$(ii)$ If $R_1>1$ and $R_2\leq1$, then $\lim\limits_{m\rightarrow\infty}\mathbf{v}(y,\ t+mT;\ \boldsymbol{\eta})=(V_1^{\vartriangle},\ M)$;

$(iii)$ If $R_1\leq1$ and $R_2>1$, then $\lim\limits_{m\rightarrow\infty}\mathbf{v}(y,\ t+mT;\ \boldsymbol{\eta})=(0,\ V_3^{\vartriangle})$;

$(iv)$ when $R_1>1$, $R_2>1$, $(\frac{a_1}{b_1})^m(1-\frac{1}{R_1})>M_2$ and $(\frac{a_2}{b_2})^m(1-\frac{1}{R_2})>M_1$, we have
\[\lim\limits_{m\rightarrow+\infty}\mathbf{v}(y,\ t+mT;\ \boldsymbol{\eta})=(\underline{V}_1,\ \underline{V}_3)(y,\ t),
\ ~ \textrm{if}~(0,\ 0)\leq\boldsymbol{\eta}\leq(\underline{V}_1,\ \underline{V}_3) ~\textrm{in}~ \Omega(0), \\[7pt]\]
and
\[\lim\limits_{m\rightarrow+\infty}\mathbf{v}(y,\ t+mT;\ \boldsymbol{\eta})=(\overline{V}_1,\ \overline{V}_3)(y,\ t),
\ ~ \textrm{if}~(\overline{V}_1,\ \overline{V}_3)\leq\boldsymbol{\eta}\leq(M_1,\ M) ~\textrm{in}~ \Omega(0).\]
\end{thm}
\bpf
$(i)$ It follows from Theorem \ref{thm2.1} that problem (\ref{1.7}) admits the unique trivial solution $(0, 0)$ when $R_1\leq1$ and $R_2\leq1$ which implies that
\begin{equation*}
  \underline{V}_1=\overline{V}_1=V_1=0
\end{equation*}
and
\begin{equation*}
  \underline{V}_2=\overline{V}_2=V_2=0.
\end{equation*}
Noticing that $V_3=M-V_2$, we have
\begin{equation*}
  \underline{V}_3=M-\overline{V}_2=M=M-\underline{V}_1=\overline{V}_3.
\end{equation*}
Recalling back to (\ref{3.5}), we have
\begin{equation*}
\begin{array}{lll}
(0,\ M)&=\liminf\limits_{m\rightarrow+\infty}\mathbf{v}(y,\ t+mT;\ \boldsymbol{\eta})\\
&\leq \limsup\limits_{m\rightarrow+\infty}\mathbf{v}(y,\ t+mT;\ \boldsymbol{\eta})=(0,\ M).
\end{array}
\end{equation*}
Thus, $\lim\limits_{m\rightarrow+\infty}\mathbf{v}(y,\ t+mT;\ \boldsymbol{\eta})$ exists and equals $(0,\ M)$.

$(ii)$ It is easy to verify that $(M_1,\ M)$ and $(0,\ M-ce^{-\lambda_2^{\vartriangle}t}\phi_2(y,\ t))$ is a pair of order upper and lower solutions of (\ref{3.2}) for some positive constant $c$ satisfying
\begin{center}
  $M-c\phi_2(y,\ 0)\leq v_{3,0}(y)\leq M$.
\end{center}
Then, it follows from Proposition \ref{prop3.1}$(i)$ that for any $\varepsilon>0$, there is a positive constant $T^*$ such that
\begin{center}
  $M-\varepsilon\leq\underline{V}_3\leq\overline{V}_3\leq M$,
\end{center}
for any $t\geq T^*$. Letting $t\rightarrow+\infty$, we have
\begin{equation*}
\left\{
\begin{array}{lll}
\underline{V}_{1t}-\frac{d_1}{\rho^2(t)}\Delta\underline{V}_1
=\underline{V}_1(a_1-\frac{n\dot{\rho}}{\rho}-c_1\underline{V}_1-b_1(M-\underline{V}_3))&\\[2mm]
\qquad \qquad \qquad \ \geq\underline{V}_1(a_1-\frac{n\dot{\rho}}{\rho}-c_1\underline{V}_1-\varepsilon),\ &y\in\Omega(0),\ t>0,\\
\underline{V}_1=0,\ &y\in\partial\Omega(0), t>0,\\
\underline{V}_1(y,\ 0)=\underline{V}_1(y,\ T)\ &y\in\Omega(0),
\end{array}
\right.
\end{equation*}
and
\begin{equation*}
\left\{
\begin{array}{lll}
\underline{V}_{1t}-\frac{d_1}{\rho^2(t)}\Delta\underline{V}_1
=\underline{V}_1(a_1-\frac{n\dot{\rho}}{\rho}-c_1\underline{V}_1-b_1(M-\underline{V}_3))&\\[2mm]
\qquad \qquad \qquad \ \leq\underline{V}_1(a_1-\frac{n\dot{\rho}}{\rho}-c_1\underline{V}_1),\ &y\in\Omega(0),\ t>0,\\
\underline{V}_1=0,\ &y\in\partial\Omega(0),\ t>0,\\
\underline{V}_1(y,\ 0)=\underline{V}_1(y,\ T)\ &y\in\Omega(0).
\end{array}
\right.
\end{equation*}
Let $\varepsilon\rightarrow0$ we have
\begin{equation}\label{3.6}
\left\{
\begin{array}{lll}
\underline{V}_{1t}-\frac{d_1}{\rho^2(t)}\Delta\underline{V}_1
=\underline{V}_1(a_1-\frac{n\dot{\rho}}{\rho}-c_1\underline{V}_1),\ &y\in\Omega(0),\ t>0,\\
\underline{V}_1=0,\ &y\in\partial\Omega(0),\ t>0,\\
\underline{V}_1(y,\ 0)=\underline{V}_1(y,\ T)\ &y\in\Omega(0).
\end{array}
\right.
\end{equation}
Similarly, we have
\begin{equation}\label{3.7}
\left\{
\begin{array}{lll}
\overline{V}_{1t}-\frac{d_1}{\rho^2(t)}\Delta\overline{V}_1
=\overline{V}_1(a_1-\frac{n\dot{\rho}}{\rho}-c_1\overline{V}_1),\ &y\in\Omega(0),\ t>0,\\
\overline{V}_1=0,\ &y\in\partial\Omega(0),\ t>0,\\
\overline{V}_1(y,\ 0)=\overline{V}_1(y,\ T)\ &y\in\Omega(0).
\end{array}
\right.
\end{equation}
According to \cite{hess}(Theorem 27.1), both (\ref{3.6}) and (\ref{3.7}) admit a unique periodic solution. Thus,  $\underline{V}_1=\overline{V}_1(\triangleq V_1^{\vartriangle})$. Recalling back to (\ref{3.5}), we have
\begin{equation*}
\begin{array}{lll}
(V_1^{\vartriangle},\ M)&=\liminf\limits_{m\rightarrow+\infty}\mathbf{v}(y,\ t+mT;\ \boldsymbol{\eta})\\
&\leq \limsup\limits_{m\rightarrow+\infty}\mathbf{v}(y,\ t+mT;\ \boldsymbol{\eta})=(V_1^{\vartriangle},\ M).
\end{array}
\end{equation*}
Thus, $\lim\limits_{m\rightarrow+\infty}\mathbf{v}(y,\ t+mT;\ \boldsymbol{\eta})$ exists and equals $(V_1^{\vartriangle},\ M)$.

$(iii)$ The proof of $(iii)$ is similar to $(ii)$.

$(iv)$ According to Theorem \ref{thm2.1}, Proposition \ref{3.1} and the  transformation $v_1=M-v_3$, we deduce that problem (\ref{3.2}) admits a minimal positive periodic solution $(\underline{V}_1,\ \underline{V}_3)$ and a maximal positive periodic solution $(\overline{V}_1,\ \overline{V}_3)$. Thus, $(M_1,\ M)$ and $(\overline{V}_1,\ \overline{V}_3)$ can be viewed as a pair of ordered upper and lower solution of (\ref{3.2}). Take
\begin{center}
  $\underline{V}_1^{(0)}=\overline{V_1},\ \underline{V}_3^{(0)}=\overline{V}_3,\ \overline{V}_1^{(0)}=M_1,\ \overline{V}_3^{(0)}=M-\delta\varphi_2$
\end{center}
as initial iterations in (\ref{3.3}). Then we have another maximal positive periodic solution of problem (\ref{3.2}) denoted by $(\overline{V}_1^{'},\ \overline{V}_3^{'})$, and another minimal positive periodic solution of problem (\ref{3.2}) denoted by $(\underline{V}_1^{'},\ \underline{V}_3^{'})$. Obviously,
\begin{center}
 $\underline{V}_1^{'}=\overline{V}_1^{'}=\overline{V}_1,\ \underline{V}_3^{'}=\overline{V}_3=\overline{V}_3^{'}$.
\end{center}
According to Proposition \ref{prop3.1}, problem (\ref{3.2}) admits the unique periodic solution $(\overline{V}_1,\ \overline{V}_3)$. And from the Lemma \ref{prop3.3}, we have
\begin{equation*}
\lim_{m\rightarrow\infty}\mathbf{v}(y,\ t+mT;\ \boldsymbol{\eta})=(\overline{V}_1,\ \overline{V}_3),
\end{equation*}
if
\begin{equation*}
  (\overline{V}_1,\ \overline{V}_3)\leq\boldsymbol{\eta}\leq(M_1,\ M).
\end{equation*}
Similarly, we have
\begin{equation*}
\lim_{m\rightarrow\infty}\mathbf{v}(y,\ t+mT;\ \boldsymbol{\eta})=(\underline{V}_1,\ \underline{V}_3),
\end{equation*}
if
\begin{center}
  $(0,\ 0)\leq\boldsymbol{\eta}\leq(\underline{V}_1,\ \underline{V}_3)$.
\end{center}
\epf

Coming back to problem (\ref{1.6}), we have the following results directly achieved from the Theorem \ref{3.4} and the  transformation $v_3=M-v_2$.
\begin{thm}\label{thm3.5}
Denote $\boldsymbol{\zeta}=(v_{1,0}(y),\ v_{2,0}(y))$ and $(v_1,\ v_2)(y,\ t;\ \boldsymbol{\zeta})$ the solution of $(\ref{1.6})$  with any nonnegative nontrivial initial value $\boldsymbol{\eta}$.

$(i)$ If $R_1\leq1$ and $R_2\leq1$, then $\lim_{m\rightarrow\infty}(v_1,\ v_2)(y,\ t+mT;\ \boldsymbol{\zeta})=(0,\ 0)$;

$(ii)$ If $R_1>1$ and $R_2\leq1$, then $\lim_{m\rightarrow\infty}(v_1,\ v_2)(y,\ t+mT;\ \boldsymbol{\zeta})=(V_1^{\vartriangle},\ 0)$;

$(iii)$ If $R_1\leq1$ and $R_2>1$, then $\lim_{m\rightarrow\infty}(v_1,\ v_2)(y,\ t+mT;\ \boldsymbol{\zeta})=(0,\ V_2^{\vartriangle})$;

$(iv)$ If $R_1>1$, $R_2>1$, $(\frac{a_1}{b_1})^m(1-\frac{1}{R_1})>M_2$ and $(\frac{a_2}{b_2})^m(1-\frac{1}{R_2})>M_1$, we have
\begin{eqnarray*}
&&\lim\limits_{m\rightarrow+\infty}(v_1,\ v_2)(y,\ t+mT;\ \boldsymbol{\zeta})\\
&=&\left\{
  \begin{array}{ll}
(\overline{v}_1,\ \underline{v}_2)(y,\ t),\ ~ \textrm{if}~(\overline{v}_1,\ 0)\leq\boldsymbol{\zeta}\leq(M_1,\ \underline{v}_2) ~\textrm{in}~ \Omega(0), \\[7pt]
(\underline{v}_1,\ \overline{v}_2)(y,\ t),\ ~ \textrm{if}~(0,\ \overline{v}_2)\leq\boldsymbol{\zeta}\leq(\underline{v}_1,\ M_2) ~\textrm{in}~ \Omega(0).
  \end{array}
  \right.
\end{eqnarray*}
\end{thm}

\section{The impact of evolution}
In order to investigate the impact of periodic evolution of domain on the competitive model, here we first present the result of (\ref{1.6}) on a fixed domain, that is (\ref{1.6}) with $\rho\equiv1$:
\begin{equation}
\left\{
\begin{array}{lll}
v_{1t}-d_1\triangle v_1=v_1(a_1-c_1v_1-b_1v_2),\ &y\in\Omega(0),\ t>0,\\
v_{2t}-d_2\triangle v_2=v_2(a_2-b_2v_1-c_2v_2),\ &y\in\Omega(0),\ t>0,\\
v_1(y,\ t)=v_2(y,\ t)=0,\ &y\in\partial\Omega(0),\ t>0,\\
v_1(y,\ 0)=v_{1,0}(y),\ v_2(y,\ 0)=v_{2,0}(y),\ &y\in\overline\Omega(0).
\end{array} \right.
\label{4.1}
\end{equation}
According to \cite{Zhaobook}(Lemma 13.1.1), the principal eigenvalue of (\ref{4.1}) is
\begin{equation}\label{4.2}
  R_i|_{\rho=1}=\frac{\int_{0}^{T}a_i(t)dt}{d_i\lambda_0\int_{0}^{T}\frac{1}{\rho^2(t)}dt}|_{\rho=1}
  =\frac{\int_{0}^{T}a_i(t)dt}{Td_i\lambda_0}(i=1,\ 2),
\end{equation}
denoted by $R_i^*$.
The corresponding periodic problem of (\ref{4.1}) is
\begin{equation}
\left\{
\begin{array}{lll}
V_{1t}-d_1\triangle V_1=V_1(a_1-c_1V_1-b_1V_2),\ &y\in\Omega(0),\ t>0,\\
V_{2t}-d_2\triangle V_2=V_2(a_2-b_2V_1-c_2V_2),\ &y\in\Omega(0),\ t>0,\\
V_1(y,\ t)=V_2(y,\ t)=0,\ &y\in\partial\Omega(0),\ t>0,\\
V_1(y,\ 0)=V_1(y,\ T),\ V_2(y,\ 0)=V_2(y,\ T),\ &y\in\overline\Omega(0).
\end{array} \right.
\label{4.3}
\end{equation}
\begin{thm}\label{thm4.1}
Denote $\bar{a}_i=\frac{1}{T}\int_0^Ta_idt(i=1,\ 2)$. There is a positive constant $D_i^*=\frac{\bar{a}_i}{\lambda_0}$ such that

$(i)$ if $d_1\in(D_1^*,\ +\infty)$ and $d_2\in(D_2^*,\ +\infty)$, $(\ref{4.3})$ admits a trivial solution which is globally asymptotically stable for problem $(\ref{4.1})$;

$(ii)$ if $d_1\in(0,\ D_1^*)$ and $d_2\in(D_2^*,\ +\infty)$, $(\ref{4.3})$ admits a semi-trivial solution, which is a global attractor for problem $(\ref{4.1})$;

$(iii)$ if $d_1\in(D_1^*,\ +\infty)$ and $d_2\in(0,\ D_2^*)$, $(\ref{4.3})$ admits a semi-trivial solution, which is global attractor for problem $(\ref{4.1})$;

$(iv)$ if $d_1\in(0,\ D_1^*)$ and $d_2\in(0,\ D_2^*)$, $(\ref{4.3})$ admits the maximal and minimal periodic solutions, which are local attractors  of problem $(\ref{4.1})$.
\end{thm}

The proof is omitted here as the assertion is easy to verified by letting $R_i^*=1$ and recalling Theorem \ref{3.5}.

Next, we consider the impact of the evolving rate on the long time behavior of the solution to problem (\ref{1.6}). There are corresponding results in the evolving domain.
\begin{thm}\label{thm4.2}
Denote $\overline{\rho^{-2}}=\frac{1}{T}\int_0^T\frac{1}{\rho^2}dt$. There is a positive constant $D_i=\frac{\bar{a}_i}{\lambda_0}\frac{1}{\overline{\rho^{-2}}}$ such that

$(i)$ if $d_1\in[D_1,\ +\infty)$ and $d_2\in[D_2,\ +\infty)$, then $(\ref{1.7})$ admits a trivial solution which is globally asymptotically stable;

$(ii)$ if $d_1\in(0,\ D_1)$ and $d_2\in(D_2,\ +\infty)$, then $(\ref{1.7})$ admits a semi-trivial solution $(V_1^{\vartriangle},\ 0)$, which is a global attractor of problem $(\ref{1.6})$;

$(iii)$ if $d_1\in(D_1,\ +\infty)$ and $d_2\in(0,\ D_2)$, then $(\ref{1.7})$ admits a semi-trivial solution  $(0,\ V_2^{\vartriangle})$, which is a global attractor of problem $(\ref{1.6})$;

$(iv)$ if $d_1\in(0,\ D_1)$ and $d_2\in(0,\ D_2)$, then $(\ref{1.7})$ admits  the maximal and minimal periodic solutions, which are local attractors of problem $(\ref{1.6})$.
\end{thm}

It can be found that $D_i$ are thresholds in terms of diffusion, and  $D_i^*$ are that in a fixed domain, and from the expressions of $D_i^*$ and $D_i$, we have the following assertions.
\begin{prop}\label{prop4.3} Recalling that $D_i^*=\frac{\bar{a}_i}{\lambda_0}$ and $D_i=\frac{\bar{a}_i}{\lambda_0}\frac{1}{\overline{\rho^{-2}}}$, we have

$(i)$ $D_i^*=D_i$ if $\overline{\rho^{-2}}=1$;

$(ii)$ $D_i^*>D_i$ if $\overline{\rho^{-2}}>1$;

$(iii)$ $D_i^*<D_i$ if $\overline{\rho^{-2}}<1$.
\end{prop}

Proposition \ref{prop4.3} implies that the evolution with a larger rate allows individuals to move with more freedom so that benefits the survival of both species, which competes each other, while the evolution with a smaller rate goes against.

\section{Numerical experiments}
In this section, Matlab is utilized to do some numerical simulations in terms of problem (\ref{1.6}) to support the theoretical results obtained in section 4. To emphasis the impact of the evolution, we assume that the diffusion rates $d_1=0.2$ and $d_2=0.1$, intrinsic population growth rates $a_1=a_2=1.2$, interspecific competition factors $b_1=b_2=0.013$, intraspecific competition factors $c_1=c_2=0.012$ and $\Omega(0)=(0,\ 1)$ followed with $\lambda _0=\pi^2$. Set the evolution rate
\begin{center}
  $\rho(t)=1-m|\sin \pi t|,\ -1<m<1$,
\end{center}
and hence
\begin{equation*}
\left\{
\begin{array}{lll}
\overline{\rho^{-2}}=1,\ &m=0,\\
\overline{\rho^{-2}}>1,\ &0<m<1,\\
\overline{\rho^{-2}}<1,\ &-1<m<0.
\end{array}\right.
\end{equation*}
Next, we select $m$ for different evolution ratios of the domain and then observe the develop trends of $v_1$ and $v_2$. The situation of $m=0$ will be presented at first for comparison.

\textbf{Example 5.1}
Set $\rho(t)=1$. Correspondingly, one has
\[\overline{\rho^{-2}}=1.\]
Meanwhile, it follows from (\ref{4.2}) that
\begin{equation*}
\left.
\begin{array}{lll}
R_1^*=\frac{\int_{0}^{T}a_1(t)dt}{Td_1\lambda_0}
&=&\frac{1.2}{0.2\pi^2}\approx0.6079<1,\\
R_2^*=\frac{\int_{0}^{T}a_2(t)dt}{Td_2\lambda_0}
&=&\frac{1.2}{0.1\pi^2}\approx1.2159>1.
\end{array}
\right.
\end{equation*}
According to Theorem \ref{thm3.4} $(iii)$, we know that $v_1$ in such fixed domain will vanish, while $v_2$ will survive. As what we have concluded, Fig. 2 (a) shows that the variable $v_2$ tends to a positive steady state while $v_1$ tends to zero, which means that the species denoted by $v_2$ will persist and $v_1$ is vanishing as time goes on.
\begin{figure}
\centering
\subfigure[]{ {
\includegraphics[width=0.30\textwidth]{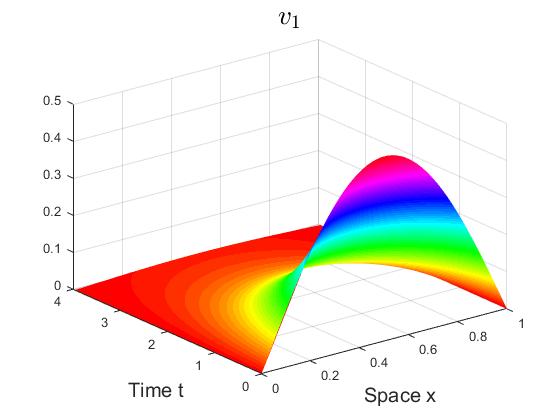}
\includegraphics[width=0.30\textwidth]{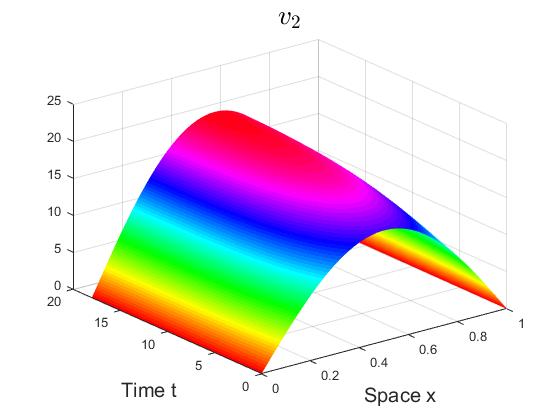}
} }
\subfigure[]{ {
\includegraphics[width=0.30\textwidth]{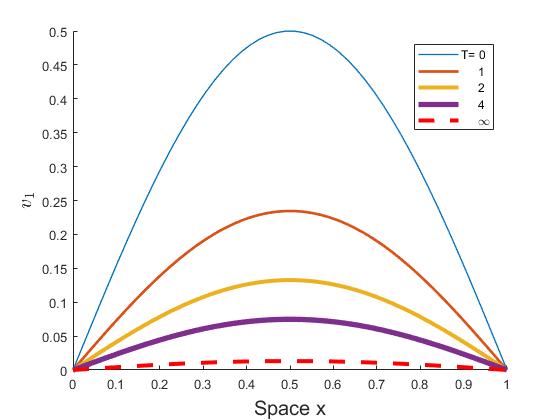}
\includegraphics[width=0.30\textwidth]{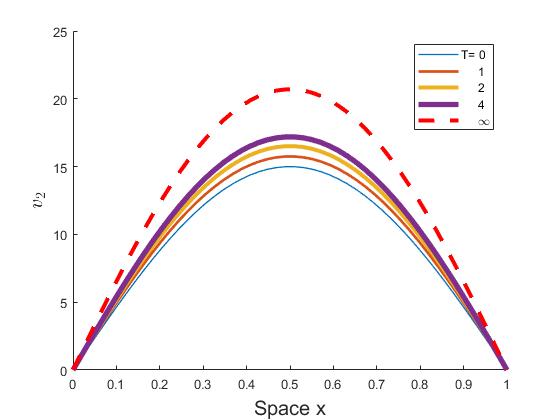}
} }
\subfigure[]{ {
\includegraphics[width=0.30\textwidth]{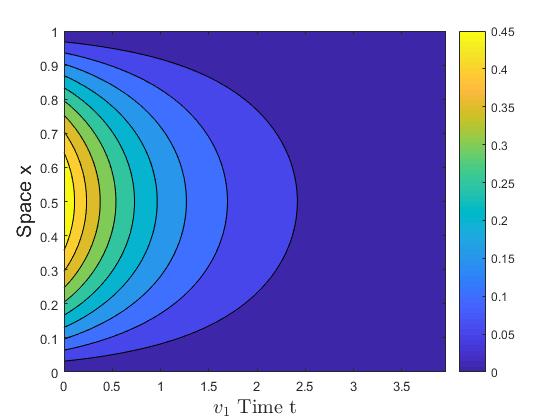}
\includegraphics[width=0.30\textwidth]{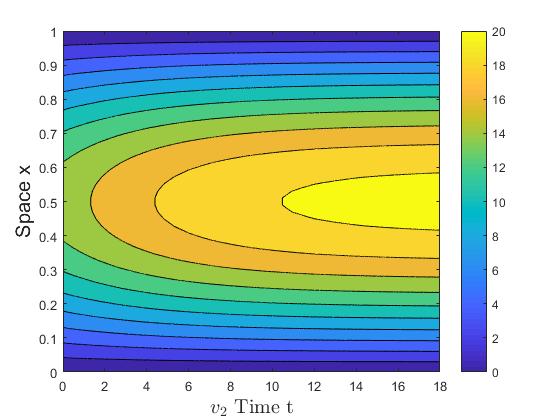}
} }
\renewcommand{\figurename}{Fig.}
\caption{\scriptsize $\rho(t)=1$. It is taken in a fixed domain. Graph (a) shows that the variable $v_2$ stabilizes to an equilibrium while $v_1$ vanishes. Graphs (b) and (c), respectively, are the cross-sectional view and contour view of graph (a). }
\end{figure}

\textbf{Example 5.2 }
Set $\rho(t)=1+0.5|\sin t|$. Correspondingly, one has
\[\overline{\rho^{-2}}=\frac{1}{2}\int_0^2\frac{1}{(1+0.5|\sin t|)^2}dt\approx0.6020.\]
Meanwhile, it follows from (\ref{2.3}) that
\begin{equation*}
\left.
\begin{array}{lll}
R_1=\frac{\int_{0}^{T}a_1(t)dt}{d_1\lambda_0\int_{0}^{T}\frac{1}{\rho^2(t)}dt}
&=&\frac{1.2}{0.2\pi^2}\frac{1}{\frac{1}{2}\int_0^2\frac{1}{(1+0.5|\sin t|)^2}dt}\\
&\approx&\frac{0.6079}{0.6020}>1,\\
R_2=\frac{\int_{0}^{T}a_2(t)dt}{d_2\lambda_0\int_{0}^{T}\frac{1}{\rho^2(t)}dt}
&=&\frac{1.2}{0.1\pi^2}\frac{1}{\frac{1}{2}\int_0^2\frac{1}{(1+0.5|\sin t|)^2}dt}\\
&\approx&\frac{1.2159}{0.6020}>1.
\end{array}
\right.
\end{equation*}
It follows from Theorem \ref{thm3.4} $(iv)$ that both $v_1$ and $v_2$ in such evolving domain will persist. As what we have concluded, Fig. 3 (a) shows that the variables $v_1$ and $v_2$ tend to  positive steady states. Fig. 3 (b) and (c) are the corresponding cross-sectional view and contour one for $v_1$ and $v_2$, respectively, and they also clearly indicate not only that the variables $v_1$ and $v_2$ keep positive, but also that the domain, to which $v_1$ and $v_2$ belong to, is periodically evolving.
\begin{figure}
\centering
\subfigure[]{ {
\includegraphics[width=0.30\textwidth]{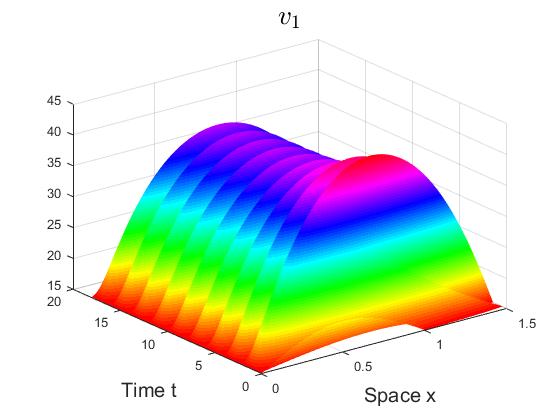}
\includegraphics[width=0.30\textwidth]{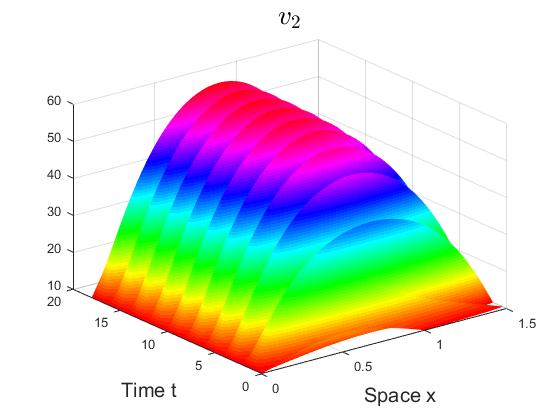}
} }
\subfigure[]{ {
\includegraphics[width=0.30\textwidth]{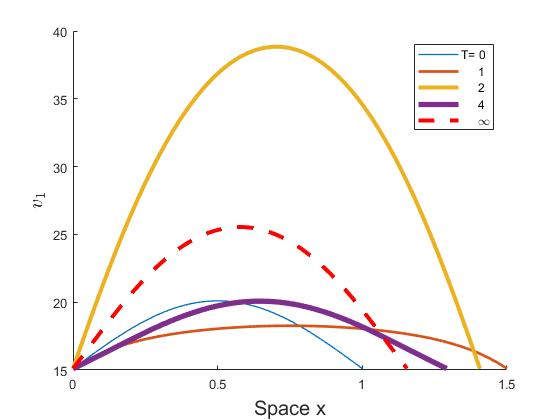}
\includegraphics[width=0.30\textwidth]{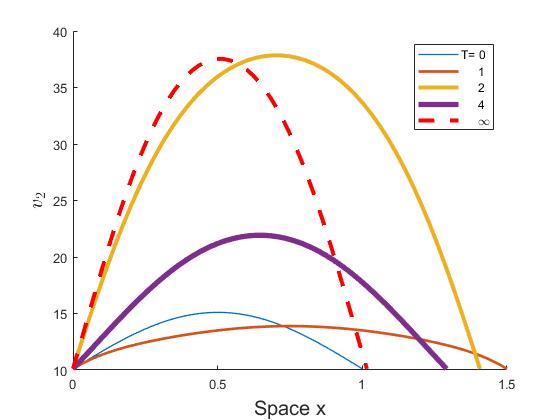}
} }
\subfigure[]{ {
\includegraphics[width=0.30\textwidth]{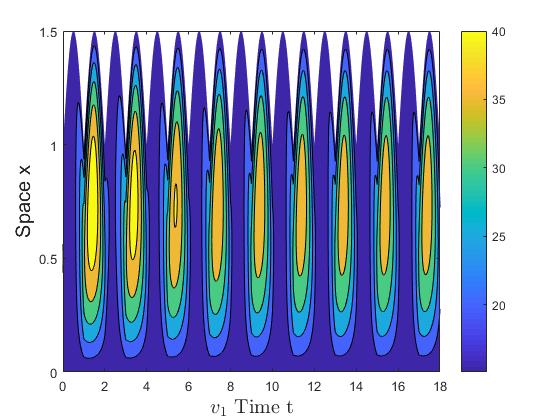}
\includegraphics[width=0.30\textwidth]{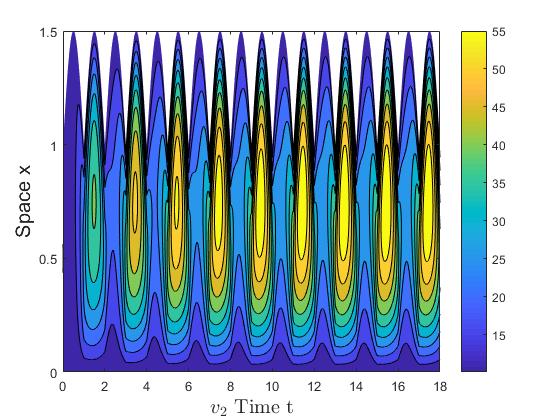}
} }
\renewcommand{\figurename}{Fig.}
\caption{\scriptsize $\rho(t)=1+0.5|\sin t|$. For the bigger evolution ratio $\rho(t)$, we acquire $R_i>1(i=1,\ 2)$, which results in the persistence of the competitive species for both sides. Graph (a) shows that both $v_1$ and $v_2$ stabilize to an equilibrium, and graphs (b) and (c) are the cross-sectional view and contour one, respectively. Also, we can clearly observe the periodic evolution of domain from (b) and (c). }
\end{figure}

\textbf{Example 5.3 }
Set $\rho(t)=1-0.3|\sin t|$. Correspondingly, one has
\[\overline{\rho^{-2}}=\frac{1}{2}\int_0^2\frac{1}{(1-0.3|\sin t|)^2}dt\approx1.5853.\]
Meanwhile, it follows from (\ref{2.3}) that
\begin{equation*}
\left.
\begin{array}{lll}
R_1=\frac{\int_{0}^{T}a_1(t)dt}{d_1\lambda_0\int_{0}^{T}\frac{1}{\rho^2(t)}dt}
&=&\frac{1.2}{0.2\pi^2}\frac{1}{\frac{1}{2}\int_0^2\frac{1}{(1-0.2|\sin t|)^2}dt}\\
&\approx&\frac{0.6079}{1.5853}<1,\\
R_2=\frac{\int_{0}^{T}a_2(t)dt}{d_2\lambda_0\int_{0}^{T}\frac{1}{\rho^2(t)}dt}
&=&\frac{1.2}{0.1\pi^2}\frac{1}{\frac{1}{2}\int_0^2\frac{1}{(1-0.2|\sin t|)^2}dt}\\
&\approx&\frac{1.2159}{1.5853}<1.
\end{array}
\right.
\end{equation*}
Similarly, Theorem \ref{thm3.4} $(i)$ tells that both $v_1$ and $v_2$ in such evolving domain will vanish. Correspondingly, Fig. 4 (a) shows that both $v_1$ and $v_2$ decay to zero which means that the species denoted by $v_1$ and $v_2$ are vanishing as time goes on, with Fig. 4 (b) and (c) reflecting the periodical evolution of the domain that $v_1$ and $v_2$ belong to.
\begin{figure}
\centering
\subfigure[]{ {
\includegraphics[width=0.30\textwidth]{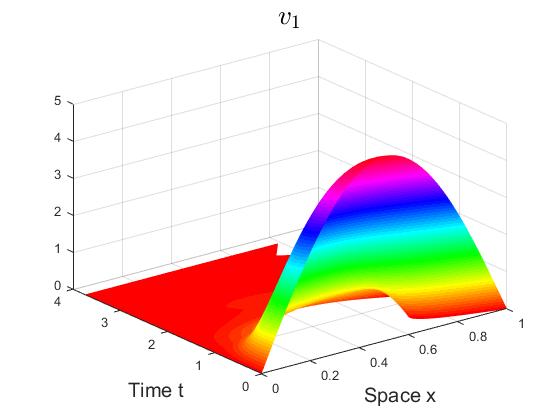}
\includegraphics[width=0.30\textwidth]{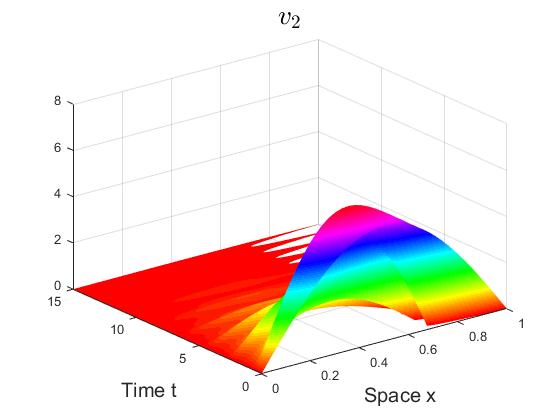}
} }
\subfigure[]{ {
\includegraphics[width=0.30\textwidth]{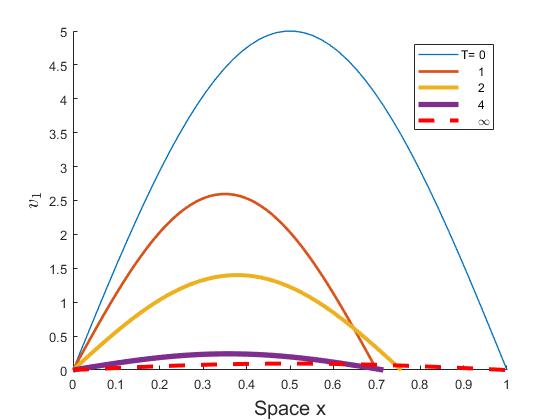}
\includegraphics[width=0.30\textwidth]{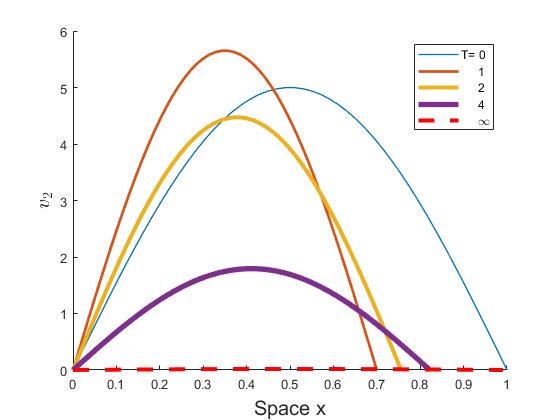}
} }
\subfigure[]{ {
\includegraphics[width=0.30\textwidth]{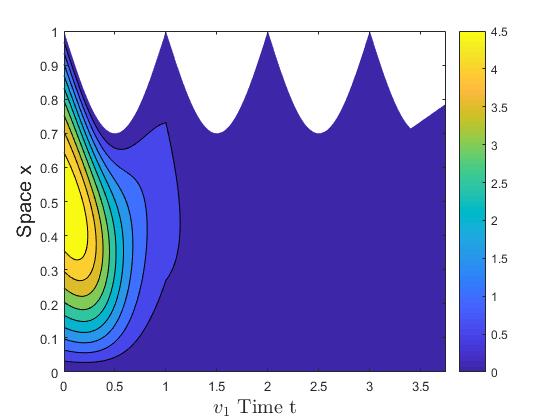}
\includegraphics[width=0.30\textwidth]{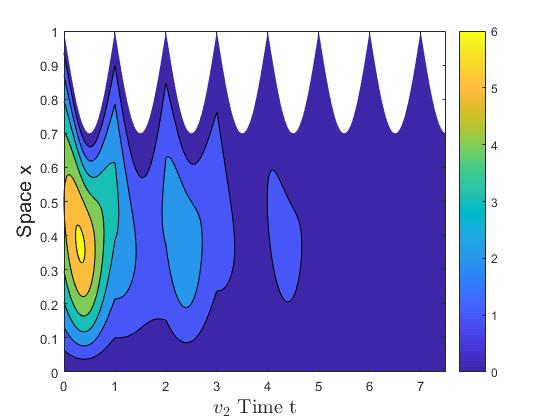}
} }
\renewcommand{\figurename}{Fig.}
\caption{\scriptsize $\rho(t)=1-0.3|\sin5t|$. For the smaller evolution ratio $\rho(t)$, we acquire $R_i<1(i=1,\ 2)$, which results in the vanishing of the two competitive species. Graph (a) shows that both $v_1$ and $v_2$ decay to the zero. Correspondingly, graphs (b) and (c) are the cross-sectional view and contour one of Graph (a), respectively, which shows the periodic evolution of the habitat. }
\end{figure}

From above, we conclude that the periodic domain evolution has a positive effect on the persistence of the species if $\overline{\rho^{-2}}<1$, but has a negative effect if $\overline{\rho^{-2}}>1$, as well as has no effect if $\overline{\rho^{-2}}=1$.

\end{document}